\documentclass[11pt]{article}
\usepackage{amsmath, amssymb, amsfonts}
\newtheorem{theorem}{Theorem}[section]

\newtheorem{proposition}[theorem]{Proposition}

\newtheorem{lemma}[theorem]{Lemma}
\newtheorem{example}[theorem]{Example}

\newtheorem{remark}[theorem]{Remark}
\numberwithin{equation}{section}

\def\IN{{\mathbb{N}}}
\def\IC{{\mathbb{C}}}
\def\IR{{\mathbb{R}}}

\def\IP{{\mathbb{P}}}

\def\b0{{\bf 0}}
\def\ba{{\bf a}}
\def\bb{{\bf b}}

\def\bx{{\bf x}}

\def\){ \right)}
\def\({\left(}
\def\Ra{{\ \Rightarrow\ }}
\def\Lra{{\ \Leftrightarrow\ }}
\def\[{\left [}
\def\]{\right ]}

\def\per{{\rm per}}
\def\qed{\hfill\vbox{\hrule width 6 pt
\hbox{\vrule height 6 pt width 6 pt}}\medskip}

\begin{document}
\openup .8\jot

\title
{Inequalities on
generalized matrix functions}

\author{Shaowu Huang\thanks{Department of Mathematics, Shanghai University, Shanghai, China.
shaowu2050@126.com}, Chi-Kwong Li\thanks{Department of Mathematics,
College of William and Mary, Williamsburg, Virginia, USA.
ckli@math.wm.edu}, Yiu-Tung Poon\thanks{Department of Mathematics,
Iowa State University, Ames, Iowa, USA. ytpoon@iastate.edu},
Qing-Wen Wang\thanks{Corresponding author. Department of Mathematics, Shanghai University,
Shanghai, China. wqw@t.shu.edu.cn     wqw369@yahoo.com}}

\date{}
\maketitle

\centerline{\bf In memory of Marvin Marcus.}

\begin{abstract} We prove inequalities on non-integer powers of products of
generalized matrices functions
on the sum of positive semi-definite matrices.
For example, for any
real number $r \in \{1\} \cup [2, \infty)$,
positive semi-definite
matrices $A_i,\ B_i,\ C_i\in M_{n_i}$, $i=1,2$, and generalized matrix functions
$d_\chi, d_\xi$ such as the determinant and permanent, etc.,
we have
\begin{eqnarray*}&&\(d_\chi(A_1+B_1+C_1)d_\xi(A_2+B_2+C_2)\)^r \\
&&\hskip 1in
 + \(d_\chi(A_1)d_\xi(A_2)\)^r + \(d_\chi(B_1)d_\xi(B_2)\)^r + \(d_\chi(C_1)d_\xi(C_2)\)^r \\
& \ge &\(d_\chi(A_1+B_1 )d_\xi(A_2+B_2 )\)^r + \(d_\chi(A_1+ C_1)d_\xi(A_2+ C_2)\)^r
+ \(d_\chi( B_1+C_1)d_\xi( B_2+C_2)\)^r\,.\end{eqnarray*}
A general scheme is introduced to prove more general inequalities
involving $m$ positive semi-definite matrices for $m \ge 3$ that
extend the results  of other authors.

\end{abstract}

AMS Classifcations. 15A15, 15A45, 15A63, 15B57.

Keywords. Positive semi-definite matrices, generalized matrix functions,
majorization.

\section{Introduction}

Let $G$ be a subgroup of the symmetric group $S_n$ of degree
$n$,  and let $\chi$ be a linear   character of
$G$. The generalized matrix function associated with $G$ and $\chi$
(also known as the $G$-immanant) of
a matrix $A=(a_{ij})\in M_n$ is defined by
$$d_{\chi}^G(A)=\sum_{\sigma\in G}\chi(\sigma)\Pi_{i=1}^n
a_{i\,\sigma(i)}.$$
For simplicity, we write $d(A) = d_\chi^G(A)$ if the
function $d_\chi^G$ is understood in the context.

Denote by $X\otimes Y$ the tensor (Kronecker) product of two
matrices $X$ and $Y$. It is known that there is a decomposable
tensor $v_\chi \in \IC^{n^2}$ such that $d(A) =  v^*_\chi (\otimes^n
A) v_\chi$; see \cite{LZ,MM,MM2}. So, one can use the theory of
tensor products and quadratic forms to study inequalities on
generalized matrix functions; see \cite{BS,Bhatia,CPZ,LS,MO,PTZ} and
their references. For example,  for any positive semi-definite
matrices $A, B \in M_n$ and any positive integer $k$, we have
$$\otimes^k (A+B) \ge \otimes^k A + \otimes^k B.$$
Letting $k = n\ell$ with $\ell \in \IN$, and using $\otimes^{\ell}
v_\chi$, we see that
\begin{equation}\label{psd2}d(A+B)^{\ell} \ge d(A)^\ell + d(B)^\ell.\end{equation}
Using similar techniques, one can obtain inequalities on
the integer powers of

\medskip
\centerline{$d(A_1), \dots, d(A_m), d(A_1+A_2), \dots, d(A_{m-1}+A_m),$ etc.}

\medskip\noindent
for positive semi-definite matrices $A_1, \dots, A_m \in M_n$.

\medskip
In this paper, we are interested in inequalities related to
non-integer powers of generalized matrix functions. In some
situations, such inequalities can be obtained by using the theory of
majorization and Schur convex functions.  Recall that for real
vectors $u, v \in \IR^n$, we say that $u$ is weakly majorized by
$v$, denoted by $u \prec_w v$, if the sum of the $k$ largest entries
of $u$ is not larger than that of $v$ for $k = 1, \dots, n$.
Furthermore, if the sums of the entries of $u$ and $v$ are the same,
we say that $u$ is majorized by $v$, denoted by $u \prec v$. A
function $f: \IR^n \rightarrow \IR$ is Schur convex if $f(u) \le
f(v)$ whenever $u \prec v$. It is known that if $f$ is Schur convex
and if $f$ is increasing in each coordinate, then $f(u) \le f(v)$
whenever $u \prec_w v$. One may see  \cite{Bhatia,MO} for the
background on majorization and Schur convex functions.

\medskip
For two positive semi-definite matrices $A$ and $B$,
we have the weak majorization relation
$$( d(A), d(B)) \prec_w (d(A+B), 0).$$
Thus, we have
$$ f(d(A+B), 0) \ge f(d(A), d(B))$$
for any Schur convex and increasing function $f$; see \cite[Chapter
II]{Bhatia} and \cite[Chatper 3]{MO}. For example, for any $p \ge
1$, the $(x_1, x_2) \mapsto x_1^p + x_2^p$ is a Schur convex
function. So, we have
$$d(A+B)^p \ge d(A)^p+ d(B)^p.$$
Actually, for any two positive semi-definite matrices $A, B \in
M_n$, it is known that (see \cite[Section IX.8.16]{Bhatia}),
$$\otimes^n (A+B)^{1/n} \ge \otimes^n A^{1/n} + \otimes^n B^{1/n}.$$
Hence,
$$d((A+B)^{1/n}) \ge d(A^{1/n}) + d(B^{1/n}).$$
As a result,  we have the weak majorization relation
$$( d(A^{1/n}), d(B^{1/n})) \prec_w (d((A+B)^{1/n} ), 0).$$
Thus, we have
$$ f(d((A+B)^{1/n})  , 0) \ge f(d(A^{1/n}), d(B^{1/n}))$$
for any Schur convex and increasing  function $f$. For example, for
any $p \ge 1$, the $(x_1, x_2) \mapsto x_1^p + x_2^p$ is a Schur
convex function. So, we have
$$d((A+B)^{1/n})^p \ge d(A^{1/n})^p+ d(B^{1/n})^p.$$
Applying this to the determinant and permanent functions, we have
$$\det(A+B)^q \ge \det(A)^q + \det(B)^q, \qquad \hbox{ for all } q \ge 1/n,$$
$$\per((A+B)^{1/n})^p  \ge \per(A^{1/n})^p + \per(B^{1/n})^p
\qquad \hbox{ for all } p \ge 1.$$
However, for three or more positive semi-definite matrices, one may
or may not be able to apply these arguments.
For example, it is known that (see \cite{BS,LS})
\begin{equation}\label{ABC}d(A+B+C) + d(A) \ge d(A+B) + d(A+C).\end{equation}
Because $(d(A+B), d(A+C)) \prec_w (d(A+B+C), d(A))$, it follows that
$$f (d(A+B+C), d(A)) \ge f(d(A+B), d(A+C))$$
for any Schur convex and increasing functions $f$. In particular,
for any $p \ge 1$,
$$d(A+B+C)^p + d(A)^p \ge d(A+B)^p + d(A+C)^p.$$
On the other hand, it is also known that
$$d(A+B+C) + d(A) + d(B) + d(C) \ge d(A+B) + d(A+C) + d(B+C).$$
However, in general,
$$(d(A+B), d(A+C), d(B+C), 0)\not\prec_w (d(A+B+C), d(A), d(B), d(C)).$$
For example, if $A = B = C = I$, then
$(d(A+B), d(A+C), d(B+C), 0) = (2^n, 2^n, 2^n, 0)$, which is not
weakly majorized by $(d(A+B+C), d(A), d(B), d(C)) = (3^n,1,1,1)$.
So, we cannot deduce that
$f(d(A+B+C), d(A), d(B), d(C)) \ge f(d(A+B), d(A+C), d(B+C), 0)$
for a Schur convex function.
We will give examples in the next section
showing that there is $p>1$ such that
\begin{equation}\label{ABC-p}
d(A+B+C)^p + d(A)^p + d(B)^p + d(C)^p \ge d(A+B)^p + d(A+C)^p + d(B+C)^p
\end{equation}
is not valid even though one can use the tensor product and quadratic form
techniques to show that (\ref{ABC-p})
holds for all positive integer $p$.

In this paper, we will develop a general scheme to prove
inequalities involving the (non-integer) powers of generalized
matrix functions. For example, we will prove in Section 2 that
(\ref{ABC-p}) is valid for any $p \in \{1\} \cup [2, \infty)$. A
general scheme and more results will be described in Section 3.
Further extensions of our techniques will be mentioned in Section 4.

\section{Results on three matrices}

Suppose we have three positive semi-definite matrices $A, B, C \in M_n$.
It was proved in \cite{BS} that
\begin{equation}
\label{oabc}\otimes^n(A+B+C) + \otimes^n A + \otimes^n B + \otimes^n C
\ge \otimes^n (A+B) + \otimes^n (A+C) + \otimes^n (B+C).\end{equation}
Applying the quadratic form using $v_\chi$ on both sides, we have
\begin{equation}
\label{abc}
d(A+B+C) + d(A) + d(B) + d(C) \ge d(A+B) + d(A+C) + d(B+C).\end{equation}
For any $\ell \in \IN$,
$$\qquad \otimes^{\ell}(\otimes^n(A+B+C)) + \otimes^{\ell}(\otimes^n A) +
\otimes^{\ell}(\otimes^n B) + \otimes^{\ell}(\otimes^n C)
$$
$$\ge \otimes^{\ell}(\otimes^n (A+B)) + \otimes^{\ell}(\otimes^n (A+C)) +
\otimes^{\ell}(\otimes^n (B+C)).$$
Applying the quadratic form using $\otimes^{\ell}(v_\chi)$ on both sides, we have
\begin{equation}
\label{integer}
d(A+B+C)^\ell + d(A)^\ell + d(B)^\ell + d(C)^\ell
\ge d(A+B)^\ell + d(A+C)^\ell + d(B+C)^\ell.
\end{equation}
As mentioned in Section 1, in general,
$$(d(A+B), d(A+C), d(B+C),0) \not\prec_w (d(A+B+C), d(A), d(B), d(C)).$$
Thus, we cannot apply the Schur convex function result to conclude that
$$d(A+B+C)^r + d(A)^r + d(B)^r + d(C)^r \ge d(A+B)^r + d(A+C)^r + d(B+C)^r$$
for $r \ge 1$.
Nevertheless,  we have the following.

\begin{theorem} \label{main}
Suppose $A, B, C\in M_n$ are positive semi-definite matrices, and $r
\in \{1\} \cup [2, \infty)$. Then for any generalized matrix
function $d(X)$, we have
\begin{equation}\label{eqr}
d(A+B+C)^r + d(A)^r + d(B)^r + d(C)^r \ge d(A+B)^r + d(A+C)^r + d(B+C)^r.
\end{equation}

\end{theorem}

\it Proof.  \rm To prove (\ref{eqr}), let
$$d(A) = x_1, \quad d(B) = x_2, \quad d(C) = x_3,$$
$$d(A+B) = x_1+ x_2 + x_{12} , \qquad
 d(A+C) = x_1 + x_3 + x_{13} ,  \qquad   d(B+C) = x_2 + x_3 +x_{23}.$$
Then $x_1, x_2, x_3 \ge 0$. Moreover,
$x_{12} = d(A+B) - \(d(A) + d(B)\) \ge 0$.
Similarly, $x_{13}, x_{23} \ge 0$.

By (\ref{abc}), we have
$$d(A+B+C) + d(A) + d(B) + d(C) \ge d(A+B) + d(B+C) + d(A+C),$$
so that
$$\begin{array}{rl}
d(A+B+C)  \ge& ( x_1+ x_2 + x_{12})+ (x_1 + x_3 + x_{13})+(x_2 + x_3 +x_{23})-(x_1 + x_2 + x_{3})\\&\\
=&( x_1+ x_2 + x_3+x_{12} +   x_{13}  +x_{23}).
\end{array}
$$
Let $X=\{( x_1, x_2 , x_3,x_{12},   x_{13} ,x_{23}):x_i,\ x_{jk}\ge 0\ \mbox{ for all }\
1\le i\le  3,\ 1\le j<k\le 3\}$. For $r\ge 2$,  define $f:X\to \IR$ by
$$f( x_1, x_2 , x_3,x_{12},   x_{13} ,x_{23})=( x_1+ x_2 + x_3+x_{12} +   x_{13}  +x_{23})^r+\sum_{i=1}^3x_i^r-\sum_{1\le j<k\le 3} (  x_j+ x_k + x_{jk})^r.$$
We have
 $$ \dfrac{\partial f}{\partial x_{12}}=r\(( x_1+ x_2 + x_3+x_{12} +   x_{13}
 +x_{23})^{r-1}-( x_1+ x_2  +x_{12} )^{r-1}\)\ge 0$$
on $X$. By symmetry,
$\dfrac{\partial f}{\partial x_{jk}}\ge 0$ on $X$ for all $1\le j<k\le 3$. Also,
 $$ \dfrac{\partial f}{\partial x_{1}}=
 r\(( x_1+ x_2 + x_3+x_{12} +   x_{13}  +x_{23})^{r-1}+x_1^{r-1} \right.$$
$$ \hskip 1in - \left. \(( x_1+ x_2  +x_{12} )^{r-1}+( x_1+ x_3  +x_{13} )^{r-1}\)\)\ge 0$$ on $X$ because

\bigskip
 $\(( x_1+ x_2  +x_{12} ),( x_1+ x_3  +x_{13} ) \)\prec_w \(( x_1+ x_2 + x_3+x_{12} +   x_{13}  +x_{23}), x_1\)$ and  $r\ge 2$.

\medskip\noindent
 By symmetry, $ \dfrac{\partial f}{\partial x_{i}}\ge 0$ on $X$ for all $1\le i\le 3$.
Therefore, for every $\bx \in X$, we have
 $$f(\bx)\ge f(\b0)=0\ \mbox{ for all }\ \bx\in X.$$
\vskip -.3in\qed

The following example shows that the bound the region for $r$  in Theorem \ref{main} is best possible for $n=1$.

\begin{example}\label{eg4} Let $A=B=C=[1]$. Then
\begin{eqnarray*}
&&\det(A+B+C)^r + \det(A)^r + \det(B)^r + \det(C)^r  \\
&& \hskip .5in   - \(\det(A+B)^r + \det(A+C)^r + \det(B+C)^r\)\\
&=&3^{r}+3-3\(2^{r}\).
\end{eqnarray*}
Let $f(r)=3^{r}+3-3\(2^{r}\) $. Then $f(1)=f(2)=0$ and
$f''(r)=3^r(\ln(3))^2 -3\(2^{r}\)(\ln (2))^2>0$ for $r\ge 1$. Therefore,  $f(r)<0$ for $1<r<2$.
\end{example}
\medskip

  The following example shows that even for $n>1$  the region for $r$  in Theorem \ref{main} cannot be extended to $[1,\ \infty)$.

\begin{example}\label{eg5} Let $A=B=\(\begin{array}{rr}1&1\\ 1&1\end{array}\)$ and $  C(x)=x\(\begin{array}{rr}1&-1\\ -1&1\end{array}\)$. Then

$$\begin{array}{c}\per(A)=\per(B)=2,\ \per(C(x))= 2x^2,\ \per(A+B)=8,\\ \\ \per(A+C(x))=  \per(B+C(x))= (1 + x)^2 + (1 - x)^2,\ \per(A+B+C)=(2 + x)^2 + (2 - x)^2.\end{array}$$
Let $x=0.17$ and $r=1.4$. Then direct calculation shows that
\begin{eqnarray*}&&
\per(A+B+C(x))^{r} + \per(A)^{r} + \per(B)^{r} +\per(C(x))^{r} \\
&& \hskip .7in
-\( \per(A+B)^{r} + \per(A+C(x))^{r} +\per(B+C(x))^{r}\) < -0.01.\end{eqnarray*}

\end{example}

\section{A general scheme and additional results}

The following observations capture the main idea  in the proof of Theorem \ref{main}.

\begin{proposition} \label{prop1} Let
$X=\{(x_1,\dots,x_N):x_i\ge 0\ \mbox{ for all }\ 1\le i\le N\}$ and
$f:X\to \IR$ be a function with continuous partial derivatives.
\begin{enumerate}
\item If
$\dfrac{\partial f}{\partial x_{j}}\ge 0$ on $X$
for all $1\le j\le N $. Then $f(\bx)\ge f(\b0)$.
\item If $(a_1, \dots, a_N), (b_1, \dots, b_N) \in X$ are such that
$(a_1, \dots, a_N) \prec_w (b_1, \dots, b_N)$ and $p \ge 1$, then
$$a_1^p + \cdots + a_N^p \le b_1^p + \cdots + b_N^p.$$
\end{enumerate}
\end{proposition}

Suppose $A_1, \dots, A_m \in M_n$ are positive semi-definite matrices
with $m \ge 3$ and $d$ a generalized matrix function on $M_n$.
For any subsequence $J$ of the sequence $K=\{1,\dots,m\}$, denoted by $J \le K$,
let $A_J=\sum_{j\in J}A_j$ and $|J|$ be the number of terms in $J$.

The following two generalizations  of (\ref{abc}) are given in
\cite[Corollary 3.4, Theorem 4.3]{BS}.

\begin{equation}\label{A1Am1}
 \sum_{j=1}^m(-1)^{m-j}\sum_{J\le K,\ |J|=j}d\(A_J\)\ge 0,
\end{equation}
\begin{equation}\label{A1Am2}
d(A_1 + \cdots + A_m)  + (m-2) \sum_{j=1}^m d(A_j)  \ge \sum_{i < j}
d(A_i+A_j).
\end{equation}
We are going to generalize (\ref{A1Am1}) and  (\ref{A1Am2}).
We continue to use the notation $I\le J$ if $I$ is a
subsequence of $J$.
We first prove the following lemma, which shows that $d(A_J)$ can be
written as the sum of non-negative numbers $x_I$ with $I \le J$.

\begin{lemma} \label{lem-3}
Suppose $d(A)$ is a generalized matrix function on $M_n$ and
$A_1, \dots, A_m \in M_n$ are positive semi-definite matrices.
Let $x_i = d(A_i) \ge 0$ for $1 \le i \le m$, and
$$x_J = d(A_J) - \sum_{L \le J,  L \ne J} x_L \qquad \hbox{ for } J \le K \hbox{ with }
|J| > 1.$$
Then $x_J \ge 0$ for every $J \le K$.
\end{lemma}

\it Proof. \rm  Let $x_J$ satisfy the hypothesis of the lemma.
We may relabel $A_1, \dots, A_m$,
and assume that $J = \{1, \dots, |J|\}$.

We prove the result by induction on $|J|$. The case for $|J|=1$ is
trivial as $d(A_1) = x_1 \ge 0$. The case for $|J|=2$ follows from
$x_{12} = d(A_1+A_2) - x_1 - x_2 = d(A_1+A_2) - d(A_1) - d(A_2) \ge
0$. Suppose the result  holds for any $q$ matrices chosen from
$\{A_1, \dots, A_{|J|}\}$ with $q < |J|$. Then we see that $x_I \ge
0$ for every $I\le J$ with $|I| < |J|$. It remains to show that $x_J
\ge 0$. By (\ref{A1Am1}),
\begin{equation} \label{a1am1-2}
\sum_{i=1}^{|J|} (-1)^{|J|-i}\sum_{I\le J,\ |I|=i}d\(A_I\)\ge
0\Ra d(A_J)\ge  -\sum_{i=1}^{|J|-1} (-1)^{|J|-i}\sum_{I\le J,\ |I|=i}d\(A_I\).
\end{equation}
Replace $d(A_I)$ by $\sum_{L \le I} x_L$ for each term on the right hand side.
After the replacement, let us determine the coefficient of $x_L$ on the right
side of (\ref{a1am1-2}) for each $L\le J$ with $1\le |L|<|J|$.
Note that $x_L$ is a summand of
$d(A_I) = \sum_{L \le I} x_L$ if and only if $L \le I$, and there are
${|J| - |L| \choose |I|-|L|}$ so many $d(A_I)$ with $L \le I$.
Thus, the  coefficient of $x_L$ on
the right side of (\ref{a1am1-2}) is
$$-\sum_{i=|L|}^{|J|-1} (-1)^{|J|-i}{|J|-|L|\choose i-|L|}
=-\sum_{j=0}^{|J|-|L|-1} (-1)^{|J|-|L|-j}{|J|-|L|\choose j} = -\((1-1)^{|J|-|L|}-1\)  = 1.$$
Therefore, we can choose
$$x_J = d(A_J) - \sum_{L \le J,  L \ne J} x_L \ge 0.$$
\qed

The following theorem extends \cite[Corollary 3.5]{BS},
which corresponds to the case when $r=1$.\medskip

\begin{theorem} \label{A1Am-1}
Let $A_1, \dots, A_m \in M_n$ be positive semi-definite matrices,
and $r \in \{1,\dots,m-2\}\cup [m-1,\ \infty)$.
Then for any generalized matrix function $d(X)$, we have
$$
 \sum_{j=1}^m(-1)^{m-j}\(\sum_{J\le K ,\ |J|=j}d\(A_J\)^r\)\ge 0,$$
 where $K=\{1,\dots,m\}$.
\end{theorem}

\it Proof. \rm  By Theorem 3.3 in \cite{BS}, for all positive integer $p$, we have
\begin{equation}\label{oA}
\sum_{j=1}^m(-1)^{m-j}\(\sum_{J\le K ,\ |J|=j}\otimes^{p} A_J  \)\ge 0\,.\end{equation}
For a positive integer $r$, let $p=nr$ and $v\in\IC^{n^2}$ such that $d(A)=v^*\(\otimes^{n}A\)v$ for all $A\in M_n$. Then the result follows by applying the quadratic form using $\otimes^rv$ on both sides of (\ref{oA}).

Suppose $r\ge m-1$.  Let  $X=\{(x_J):x_J\ge 0,\   J\le K\}$,
$K_{j}=\{J\le K:|J|=j\}$ for $1\le j\le m$,  and
$s_J = \sum_{I \le J} x_I$ for $J\in K_{j}$. Suppose
\begin{equation}\label{fbx}f(\bx)=\sum_{j=1}^m(-1)^{m-j}\sum_{J\in
K_j}s_J^r.
\end{equation}
We will prove that $\dfrac{\partial f}{\partial x_I}\ge 0$ for all
$I\le  K$, so that $f(\bx)\ge f(\b0)=0$ for $\bx\in X$. It turns out
that we need to use the nonnegativity of the higher order partial
derivatives of $f$ to prove that the first partial derivatives of
$f$ are nonnegative, which corresponds to  the case when $t =1$ in
the following claim.

\medskip\noindent
{\bf Claim: } Let $1\le t\le m-1$ and $I_1,\dots, I_t\le K$ such
that $I_j\setminus \cup_{i=1}^{j-1}I_i\neq \emptyset$ for $2\le j\le
t$. Then
$$ \dfrac{\partial^t f}{\partial x_{I_1}\dots\partial x_{I_t} }\ge 0.$$

Note that for $I,\ J\le K$, $\dfrac{\partial s_J^r}{\partial x_I} = rs_J^{r-1}$ if $I\le J$ and $0$ otherwise. More generally, for $I_1,\ I_2,\dots,\ I_t,\ J\le K$, $$\dfrac{\partial^t s_J^r}{\partial x_{I_1}\dots\partial x_{I_t} } = r(r-1)\dots (r-t+1)s_J^{r-t}\qquad \mbox{ if }\cup_{j=1}^tI_j\le J\,,$$
and $0$ otherwise.

For $I\le K$, let $X_I=\{\bx\in X: x_J=0\ \mbox{for all }J\not \le I\}$. Then for every $\bx\in X_I$ and $I\le J\le K$, we have $s_J(\bx)=s_{I}(\bx)$.
Let $I=\cup_{i=1}^{t}I_i$ and $\bx\in X_I$. Then if $\cup_{j=1}^tI_j\le J$, we have
$$ \dfrac{\partial^t s_J^r}{\partial x_{I_1}\dots\partial x_{I_t} }(\bx)=r(r-1)\dots (r-t+1)s_J^{r-t}(\bx)=r(r-1)\dots (r-t+1)s_I^{r-t}(\bx)\,.$$
Let $|I|=p$. We have
\begin{equation}\label{dxJ}
\begin{array}{rl}& \dfrac{\partial^t f}{\partial x_{I_1}\dots\partial x_{I_t} }(\bx)\\&\\
=&r(r-1)\dots (r-t+1)\(\sum_{j=p}^m(-1)^{m-j}\sum_{J\in K_j,\ I\le J}s_J^{r-t}\)(\bx)\\&\\
=&r(r-1)\dots (r-t+1)\(\sum_{j=p}^m(-1)^{m-j}\sum_{J\in K_j,\ I\le J}s_I^{r-t}\)(\bx)\\&\\
=&r(r-1)\dots (r-t+1)\(\sum_{j=p}^m(-1)^{m-j}\(\begin{array}{c} m-p\\ j-p\end{array}\)\)
s_I^{r-t}(\bx)\\&\\
=&(-1)^{m-p}r(r-1)\dots
(r-t+1)\(\sum_{j=0}^{m-p}(-1)^{j}\(\begin{array}{c} m-p\\
j\end{array}\)\)
s_I^{r-t}(\bx)\\&\\
=&0.\end{array}
\end{equation}
We will prove the claim by (backward) induction on $t$. For $t=m-1$,
by the condition on $I_i$, we have $|I|=m-1$ or $m$. If $|I|=m-1$,
then
$$ \dfrac{\partial^t f}{\partial x_{I_1}\dots\partial x_{I_t} }
=r(r-1)\dots (r-m+2)\(s_K^{r-m+1}-s_I^{r-m+1}\)\ge 0\,$$
because $r\ge m-1$.
If $|I|=m$, then
$$ \dfrac{\partial^t f}{\partial x_{I_1}\dots\partial x_{I_t} }=r(r-1)\dots (r-m+2) s_K^{r-m+1} \ge 0\,.$$
Suppose the result holds for some $1<t\le m-1$. Let $I_1,\dots I_{t-1}\le K$ such that $I_j\setminus \cup_{i=1}^{j-1}I_i\neq \emptyset$ for $2\le j\le t-1$. Let $I=
\cup_{i=1}^{t-1}I_i$. For each $I_t\le K$, with $I_t\not\le I$, we have

\begin{equation}\label{It}\dfrac{\partial \ }{\partial x_{I_t}} \( \dfrac{\partial^{t-1} f}{\partial x_{I_1}\dots\partial x_{I_{t-1}} } \)=\dfrac{\partial^t f}{\partial x_{I_1}\dots\partial x_{I_t} }\ge 0.\end{equation}

For $\bx=(x_J)\in X$, let $\bx_I=(x'_{J})\in X_I$, where $x'_J=x_J$ if $J\le I$ and $0$ otherwise. By (\ref{It}), we have

$$\dfrac{\partial^{t-1} f}{\partial x_{I_1}\dots\partial x_{I_{t-1}} }(\bx)\ge \dfrac{\partial^{t-1} f}{\partial x_{I_1}\dots\partial x_{I_{t-1}} }(\bx_I)=0$$
by (\ref{dxJ}).
\qed

\begin{example} Suppose  $m\ge 3$ and  $A_i=[1]$ for $1\le i\le m$. Let $r>0$

$$
f(m,r)= \sum_{j=1}^m(-1)^{m-j}\(\sum_{J\le K ,\ |J|=j}d\(A_J\)^r\)
= \sum_{j=1}^m(-1)^{m-j}{m\choose j}j^r.$$
Then $f(m,r)$ is the $m^{\rm th}$ finite difference (with step size = 1) of the function
$g(x)=x^r$ at $0$. By the mean value theorem of finite difference \cite{SR}, $f(m,r)$ has
the same sign as $g^{(m)}(x)$ for $x>0$. Therefore,
$f(m,r)=0$ for $r=1,2,\dots,m-1$ and  $f(m,r)<0$ for $m-2i<r<m-2i+1$ with $1\le i\le [m/2]$ and
$f(m,r)>0$ for $m-2i-1<r<m-2i$ with $1\le i\le [(m-1)/2]$. Hence, the condition on $r\ge m-1$
in Theorem \ref{A1Am-1} is necessary.
\end{example}

The next result extend the inequality in  \cite[Theorem 4.8]{BS} to non-integer powers.

 \begin{theorem} \label{A1Am-4}
Suppose $A_1, \dots, A_m \in M_n$ are positive semi-definite matrices,
and $r \in \{1\}\cup [2,\ \infty)$.
Let $K=\{1,\dots,m\}$ and $K_j=\{J\le K:|J|=j\}$ for $1\le j\le m$. For each $J\le K$,
let $A_J=\sum_{j\in J}A_j$.
Then for every $r\ge 2$, $1\le k<\ell< p\le m$ and any generalized matrix function $d(A)$,
we have
\begin{equation}\label{pkl} \dfrac{\ell-k}{p\(\begin{array}{c} m\\ p\end{array}\)}
\sum_{J\in K_p} d\( A_J\)^r +
\dfrac{p-\ell}{k\(\begin{array}{c} m\\ k\end{array}\)}
\sum_{J\in K_k} d\(A_J\)^r
\ge\dfrac{p-k}{\ell\(\begin{array}{c} m\\
\ell\end{array}\)} \sum_{J\in K_\ell} d\( A_J \)^r.\end{equation}
\end{theorem}

To prove the theorem, one needs only consider the case when $m = p$
for the following reason.
If the special case when $m = p$ is proved, then for any
$\hat J \in K_p$, we have
\begin{equation}\label{hatJ}
\dfrac{\ell-k}{p}
d\( A_{\hat J} \)^r +
\dfrac{p-\ell}{k\(\begin{array}{c} p\\ k\end{array}\)}
\sum_{J\in K_k, J \le \hat J} d\(A_J\)^r
\ge\dfrac{p-k}{\ell\(\begin{array}{c} p\\
\ell\end{array}\)} \sum_{J\in K_\ell, J \le \hat J} d\( A_J \)^r.
\end{equation}
Note that every $J \in K_k$ will appear
in ${m-k \choose p-k}$ different $\hat J \in K_p$, and
every $J \in K_\ell$ will appear in ${m-\ell \choose p-\ell}$ many
$\hat J \in J_p$.
Hence, summing the inequalities (\ref{hatJ}) for different
$\hat J$, and dividing the resulting inequality
by ${m\choose p}$, we have
 $$
  \dfrac{\ell-k}{p{m\choose p}}
\sum_{J\in K_p} d\( A_J\)^r +
\dfrac{(p-\ell){m-k\choose p-k}}{k{m\choose p}{p\choose k}}
\sum_{J\in K_k} d\(A_J\)^r
\ge\dfrac{(p-k){m-\ell \choose p-\ell}}{\ell{m\choose p}{p\choose \ell}}
\sum_{J\in K_\ell} d\( A_J \)^r,
$$
that simplifies to (\ref{pkl}).

To prove Theorem \ref{A1Am-4}, we first
consider the special case where $\ell =p-1$ and $k=p-2$ in the following lemma.

\begin{lemma}\label{pm}
Let $X=\{(x_J):x_J\ge 0,\   J\le K\}$. For $3\le p\le m$, define $f$
on $X$ by
 $$f(\bx)=(m-p)!(p-1)!\sum_{J\in K_{p}}s_J^r+(m-p+2)!(p-3)!\sum_{J\in K_{p-2}} s_J^r-2(m-p+1)!(p-2)!\sum_{J\in K_{p-1}} s_J^r\,.$$
Then $f(\bx)\ge 0$ for all $\bx\in X$.
\end{lemma}

\it Proof.
\rm By the comment after Theorem \ref{A1Am-4}, it
suffices to prove the case when $p=m$. We need to show that for all $\bx\in X$,
 $$\begin{array}{rl}&(m- 1)!s_K^r+
 2 !(m-3)!\sum_{J\in K_{m-2}} s_J^r-2(m- 2)!\sum_{J\in K_{m-1}} s_J^r\ge 0\\&\\
 \Lra&g(\bx)=(m- 1)(m-2)s_K^r+ 2  \sum_{J\in K_{m-2}} s_J^r-2(m- 2) \sum_{J\in K_{m-1}} s_J^r\ge 0.\end{array}
 $$
 For $1\le q\le m$, let $K_q=\{J\le K:|J|=q\}$.  Let $I\in K_q$.  Then
$$\dfrac{\partial g}{\partial x_I} =r\((m- 1)(m-2)s_K^{r-1}+ 2  \sum_{J\in K_{m-2},\ I\le J} s_J^{r-1}-2(m- 2) \sum_{J\in K_{m-1},\ I\le J} s_J^{r-1}\).$$
There are $(m-q)$ $J\in K_{m-1}$ such that $ I\le J$. For each such $J$,
we have $s_K\ge s_J$. If \newline$(m- 1)(m-2)\ge 2(m- 2)(m-q)\Lra q\ge \dfrac{m+1}{2}$,
then we have $\dfrac{\partial g}{\partial x_I}\ge 0$. Let
$$X_0=\{\bx\in X: x_I=0\
\mbox{ for all } I\mbox{ with  }|I|\ge \dfrac{m+1}{2}\}\,.$$ For $\bx\in X$, let $\bx_0$
be the projection of $\bx $ to $X_0$. Then we have $g(\bx)\ge g(\bx_0)$. It suffices to
prove that
\begin{equation}\label{gI}
\dfrac{\partial g}{\partial x_I}(\bx_0)\ge 0\qquad \mbox{ for all }\bx_0\in X_0,\ \mbox{ and }I\le K  \mbox{ with }|I|<\dfrac{m+1}{2}\,.
\end{equation}
Then we  have $g(\bx)\ge g(\bx_0)\ge g(\b0)=0$.

Note that for $\bx\in X_0$, $s_J(\bx)= \sum_{  I\le J,\
|I|<\frac{m+1}{2}} x_I$. Without loss of generality, we may assume
that  $I=\{m-q+1,\dots,m\}$ with $q<\dfrac{m+1}{2}$.

For $1\le i  \le m-q$, ($1\le i< j\le m-q$),  let $K(i)=K\setminus
\{ i\}$ ($K(i,j)=K\setminus \{i,\ j\}$). Let $s(i)=s_{K(i)}$ and
$s(i,j)=s_{K(i,j)}$. We may assume that $s(1)\ge s(2)\ge\cdots\ge
s(m-q)$. Let $\ba=(a_i)$ and $\bb=(b_i)\in \IR^{M}$, where
$M=(m-1)(m-2)+2{m-q\choose 2}=(m-1)(m-2)+(m-q)(m-q-1)$ be given by
$$\begin{array}{rl}
\ba&= (\underbrace{s_K,\dots ,s_K}_{(m-1)(m-2) \mbox{ \tiny terms }},
s(1,2),s(1,2),s(1,3),s(1,3),s(2,3),s(2,3),\cdots \\
& \hskip 2in \cdots,  s(1,m-q), \dots  ,s(m-q-1,m-q) ),\\&\\
\bb&= (\underbrace{s(1),\dots ,s(1)}_{2(m-2) \mbox{ \tiny terms }},
\underbrace{s(2),\dots ,s(2)}_{2(m-2) \mbox{ \tiny terms }},\dots,
\underbrace{s(m-q),\dots ,s(m-q)}_{2(m-2) \mbox{ \tiny terms }} ,\hskip-.2in
\underbrace{0,\dots,0}_{ (q-1)(q-2) \mbox{ \tiny terms
}}\hskip-.2in).\end{array}$$ Note: $M-2(m-2)(m-q) =(q-1)(q-2)$.

We show that $\bb\prec_w\ba$ in the following.
For each $1\le k\le M$, $\sum_{i=1}^ka_i=\sum_{J\le K} n(a,\,J)\,x_J$ and  $\sum_{i=1}^kb_i=\sum_{J\le K} n(b,\,J)\,x_J$ for some non-negative integers $n(a,\,J),\ n(b,\,J)$. Since $x_J\ge 0$ for all $J\le K$, it suffices to show that $n(a,\,J)\ge n(b,\,J)$ for all $J\le K$.

For each $J\le K$ with $|J|<\dfrac{m+1}{2}$,  note that every $s_K$ contains a copy of $x_J$ and  each $s(i)$ (respectively, $s(i,j)$) contains a copy of $s_J$ if and only if $i\not \in J$  (respectively, both $i,\ j\not\in J$).
Consider the following cases:

\medskip
{\bf Case 1: } $1\le k\le (m-1)(m-2)$. In this case, $n(a,\,J)=k\ge n(b,\,J)$.

\medskip
{\bf Case 2: } $(m-1)(m-2)<k\le M$. We may assume that $k\le  2(m-2)(m-q)$. Choose $t_1>1 $ such that $(t_1-2)(t_1-1)<k-(m-1)(m-2)\le (t_1-1)t_1$ and $t_2\ge 1 $ such that $2(m-2)(t_2-1)<k\le 2(m-2)t_2$. Let $k_1=k-(m-1)(m-2)-(t_1-2)(t_1-1)$.  Then $ \sum_{i=1}^ka_i$ consists of $(m-1)(m-2)$ copies of $s_K$ and 2 copies each of $s(i,j)$, where $1\le i<j\le t_1-1$ and 2 copies each of $s(i,t_1)$, where $1\le i\le \[\dfrac{k_1}{2}\]$ and a copy of  $s(\dfrac{k_1+1}{2},t_1)$ if $k_1$ is odd. On the other hand, $  \sum_{i=1}^kb_i$ consists of $2(m-2)$ copies of $s(i)$ for $1\le i\le t_2-1$ and $k- 2(m-2)(t_2-1)$ copies  of $s(t_2)$.

For every integer $m,\ t$, we have
\begin{equation}\label{mt}(m-1)(m-2)+t(t-1) = 2(m-2)t+(m-(t+1))(m-(t+2))\ge 2(m-2)t\,.\end{equation}
In particular, we have

$$\begin{array}{rl}
 &(m-1)(m-2)+t_2(t_2-1)\ge 2(m-2)t_2\\&\\
 \Ra&t_2(t_2-1)-(k-(m-1)(m-2))\ge 2(m-2)t_2-k\ge 0\\&\\
 \Ra&t_2(t_2-1)\ge  k-(m-1)(m-2)\,. \end{array}$$
 This shows that $t_1\le t_2$.

\medskip
 If $J\le I$, then $n(a,\,J)=k=n(b,\,J)$.

\medskip
 Suppose $J\not\le I$. Let $J\setminus I=\{j_1,\dots,j_u\}$ with $j_1<j_2<\cdots <j_v\le t_2-1<t_2\le j_{v+1}<\cdots<j_u$.
 Then
 $$n(b,\,J)=\left\{\begin{array}{ll}k-2(m-2)v&\mbox{ if }j_{v+1}\neq t_2\\&\\
 2(m-2)(t_2-1-v)&\mbox{ if }j_{v+1}=t_2.\end{array}\right.$$
 On the other hand,
 $$n(a,\,J)= k- \mbox{ number of copies of } s(i,j) \mbox{ in }\sum_{i=1}^ka_i \mbox{ such that either  }  i\mbox{ or  } j\in J\setminus I.$$
Consider the following subcases:

\noindent
 {\bf Case 2a: } $t_1<t_2$ or $t_2< j_{v+1}$. Then  $ t_1< j_{v+1}$. We have
 $$\begin{array}{rl}n(a,\,J)
 \ge &k-2v(t_1-1)\\&\\
 \ge &k-2v(m-q-1)\\&\\
 \ge&k-2v(m-2)\\&\\
 \ge&n(b,\,J).
 \end{array}
 $$
{\bf Case 2b: } $t_1=t_2$ and $t_2= j_{v+1}$. Then by (\ref{mt}),

$$k_1=k-(m-1)(m-2)-(t_1-2)(t_1-1)\le k-2(m-2)(t_1-1) =k-2(m-2)(t_2-1).$$

So we have

 $$\begin{array}{rl}n(a,\,J)
 \ge &k-2v(t_1-2)-k_1\\&\\
 \ge &k-2v(m-2)- (k-2(m-2)(t_2-1))\\&\\
 =& 2(m-2)(t_2-1-v)\\&\\
=&n(b,\,J).
 \end{array}
 $$
\vskip -.3in \qed

\medskip
\it Proof of Theorem \ref{A1Am-4}. \rm For $1\le q\le m$, let
$t_q=\dfrac{1}{q{m\choose q}}
\sum_{J\in K_q} d(A_{J})^r$. Then (\ref{pkl}) is equivalent to

\begin{equation}\label{tpkl}
(\ell-k)(t_p-t_\ell)\ge (p-\ell)(t_\ell-t_k)
\end{equation}
for all $1\le k<\ell<p\le m$.

By Lemma \ref{pm}, (\ref{tpkl}) holds for $\ell =p-1$ and $k=p-2$. Thus, for every $p\ge q>\ell\ge q'>k$, we have
$$t_q-t_{q-1}\ge t_{q-1}-t_{q-2}\ge \cdots \ge t_{q'}-t_{q'-1}.$$
 Therefore, we have
 $$\begin{array}{rl}
 &(\ell-k)(t_q-t_{q-1})\ge
 \sum_{q'=k+1}^\ell (t_{q'}-t_{q'-1})=t_{\ell}-t_k\\&\\
 \Ra&(\ell-k)(t_p-t_\ell)=(\ell-k)\sum_{q=\ell+1}^p(t_q-t_{q-1})\ge(p-\ell)(t_{\ell}-t_k).
 \end{array}
 $$

\vskip -.3in
\qed
\begin{remark} When $m=p=3,\ \ell=2,\ k=1$, (\ref{pkl}) reduces to (\ref{eqr}). Therefore, Example \ref{eg4} shows that the condition $r\ge 2$ in Theorem \ref{A1Am-4} is best possible.
\end{remark}

 \section{Additional results and techniques}

For any partition $\{I_1, \dots, I_k\}$ of
$\{1, \dots, m\}$, and for any Schur Convex function $f: \IR^k \rightarrow \IR$,
we have
 $$f(d(A_1+ \cdots + A_m), 0, \dots, 0))
 \ge f(d(A_{I_1}, \cdots ,d(A_{I_k}))).$$
More generally, suppose $(I_1, \dots, I_k)$ and
$(J_1, \dots, J_k)$ are two collections of muti-subsets
$\{1, \dots, m\}$, we can define
$(I_1, \dots, I_k) \triangleleft (J_1, \dots, J_k)$ if
the union of $\ell$ subsets in the family $(I_1, \dots, I_k)$
is always contained in the union of $\ell$ subsets in
the family $(J_1, \dots, J_k)$ for $\ell = 1, \dots, k$.
We have the following.

\begin{theorem}
Let $A_1, \dots, A_m \in M_m$ be positive semi-definite matrices.
Suppose $(I_1, \dots, I_k)$ and
$(J_1, \dots, J_k)$ are two collections of multi-subsets
$\{1, \dots, m\}$ such that
$(I_1, \dots, I_k) \triangleleft (J_1, \dots, J_k)$.
Then for any Schur Convex function $f: \IR^k \rightarrow \IR$,
we have
 $$f(d(A_{J_1}) + \cdots + d(A_{J_k}))
 \ge f(d(A_{I_1}, \cdots d(A_{I_k})).$$
\end{theorem}

One can also take partial sum of the positive semi-definite
matrices $A_1, \dots, A_m$ in Theorem \ref{A1Am-4}, and obtain the following
result that
removes the restriction   $r\ge 2$.

 \begin{theorem} \label{A1Am-5}
Suppose $A_1, \dots, A_m \in M_n$ are positive semi-definite matrices,
and $\Phi $ is a convex function on $[0,\ \infty)$.
Let $K=\{1,\dots,m\}$ and $K_j=\{J\le K:|J|=j\}$ for $1\le j\le m$.
Then for every $1\le k<\ell< p\le m$ and any generalized matrix function $d(X)$,
we have
\begin{equation}\label{pkl_1} \dfrac{\ell-k}{ \(\begin{array}{c} m\\ p\end{array}\)}
\sum_{J\in K_p} \Phi\(d\( A_J\)\) +
\dfrac{p-\ell}{ \(\begin{array}{c} m\\ k\end{array}\)}
\sum_{J\in K_k} \Phi\(d\(A_J\)\)
\ge\dfrac{p-k}{ \(\begin{array}{c} m\\
\ell\end{array}\)} \sum_{J\in K_\ell} \Phi\(d\( A_J \)\).\end{equation}
\end{theorem}

\it Proof. \rm Suppose $J\in K_p$. Then for every $I\in K_{p-2} $ with $I\le J$, let $J\setminus I=  \{j_1,\ j_2\}$ and  $J_i=I\cup\{j_i\}$ for $i=1,2$. Then $J_1,\ J_2$ are the only  $\hat J\in K_{p-1}$ such that $I\le \hat J\le J$. By (\ref{ABC}), we have $\(d\(A_{J_1}\),d\(A_{J_2}\)\)\prec_w\(d\(A_{J}\),d\(A_{I}\)\)$, we have $\Phi\(d\(A_{J}\)\)+\Phi\(d\(A_{I}\)\)\ge \Phi\(d\(A_{J_1}\)\)+\Phi\(d\(A_{J_2}\)\)$. Summing over all $I\in K_{p-2}$, with $I\le J$ we have
  $${p\choose p-2}\Phi\(d\(A_{J}\)\)+\sum_{I\in K_{p-2},\ I\le J}\Phi\(d\(A_{I}\)\)\ge \dfrac{2{p\choose p-2}}{p}\sum_{\hat J\in K_{p-1},\  \hat J\le J} \Phi\(d\(A_{\hat J}\)\).$$
Then summing over all $J\in K_p$, we have
\begin{equation}\label{p2}
\begin{array}{rl}
&{p\choose p-2}\sum_{J\in K_p}\Phi\(d\(A_{J}\)\)+{m-(p-2)\choose  2}\sum_{I\in K_{p-2}}\Phi\(d\(A_{I}\)\)\\&\\
\ge& (p-1)(m-(p-1))\sum_{\hat J\in K_{p-1} } \Phi\(d\(A_{\hat
J}\)\).\end{array}
\end{equation}
For $1\le j\le m$, let  $t_j=\dfrac{1}{{m\choose j}}\sum_{J\in K_j}\Phi\(d\(A_{J}\)\)$, then (\ref{p2}) is equivalent to
$$t_p+t_{p-2}\ge 2t_{p-1}\Lra t_p-t_{p-1}\ge t_{p-1}-t_{p-2}\,.$$
Following the argument at the end of the proof of Theorem \ref{A1Am-4},
we have
$$(\ell-k)(t_p-t_\ell) \ge(p-\ell)(t_{\ell}-t_k),$$
which is equivalent to (\ref{pkl_1}).
\qed

Let $\IP_k$ be the set of functions $\Phi$ on $[0,\ \infty)$ such that $\Phi^{(i)}(x)\ge 0$ for all $0\le i\le k$ and $x\ge 0$. Then the term $d(A_J)^r$ in Theorem \ref{main} (respectively, Theorem \ref{A1Am-1} and Theorem \ref{A1Am-4}) can be replaced by $\Phi(d(A_J))$ for all $\Phi\in \IP_2$ (respectively, $\IP_{m-1}$ and $\IP_2$).

\bigskip
Finally, we point out that following the same proof in \cite{BS}, (\ref{oabc}) can be generalized to the following:

\begin{proposition}
Suppose $A_i, B_i, C_i\in M_{n_i}$ are positive semi-definite matrices for $1\le i\le k$. Then
$$\otimes_{i=1}^k(A_i+B_i+C_i) + \otimes_{i=1}^k A_i + \otimes_{i=1}^k B_i + \otimes_{i=1}^k C_i
$$
$$
\ge \ \otimes_{i=1}^k (A_i+B_i) + \otimes_{i=1}^k (A_i+C_i) + \otimes_{i=1}^k (B_i+C_i). \ \ \
$$
Consequently, we have
\begin{equation}
\label{oabci}\begin{array}{rl}&\otimes_{i=1}^k\(\otimes^{n_i}(A_i+B_i+C_i)\) + \otimes_{i=1}^k\(\otimes^{n_i} A_i\) + \otimes_{i=1}^k\(\otimes^{n_i} B_i\) + \otimes_{i=1}^k\(\otimes^{n_i} C_i\)\\&\\
\ge& \otimes_{i=1}^k\(\otimes^{n_i} (A_i+B_i)\) +
\otimes_{i=1}^k\(\otimes^{n_i} (A_i+C_i)\) + \otimes_{i=1}^k
\(\otimes^{n_i}(B_i+C_i)\).\end{array}\end{equation}
\end{proposition}

For given generalized matrix functions $d_i$ on $M_{n_i}$, we can choose unit vectors
$v_i\in \IC^{n_i^2}$ such that $d_i(X)=v_i^*Xv_i$ for $X\in M_{n_i}$.
Let $v=\otimes_{i=1}^kv_i$ and apply the quadratic form using  $v$ on both
sides of (\ref{oabci}).  We have
$$\begin{array}{rl}&\otimes_{i=1}^kd_i(A_i+B_i+C_i)  + \otimes_{i=1}^kd_i(A_i)  +
\otimes_{i=1}^kd_i(B_i)  + \otimes_{i=1}^kd_i(C_i) \\&\\
 \ge &\otimes_{i=1}^kd_i(A_i+B_i)  + \otimes_{i=1}^kd_i(A_i+C_i)  +
 \otimes_{i=1}^kd_i(B_i+C_i).\end{array} $$
Thus we can replace $d(X)$ in  Theorem \ref{main}
 by a product of generalized matrix functions
$d(X_1,\dots,X_k)= d_1(X_1) \cdots d_k(X_k)$.
For example, if we set  $d(X_1,X_2) = \det(X_1)\per(X_2),$  then for any
positive semi-definite matrices $A_1, B_1, C_1\in M_{n_1}$,  $A_2, B_2, C_2\in M_{n_2}$ and $r \in \{1\} \cup [2, \infty)$,
we have
$$\begin{array}{rl}& [\det(A_1+B_1+ C_1)\per(A_1+B_1+ C_1)]^r \\
& \\
&\hskip 1in  + [\det(A_1)\per(A_2)]^r + [\det(B_1)\per(B_2)]^r +[\det(C_1)\per(C_2)]^r \\&\\
 \ge&
[\det(A_1+B_1)\per(A_2+B_2)]^r + [\det(A_1+C_1)\per(A_2+C_2)]^r +
[\det(B_1+C_1)\per(B_2+C_2)]^r.\end{array}$$

Theorem \ref{A1Am-1}, Theorem \ref{A1Am-4} and Theorem \ref{A1Am-5} can also
be generalized in a similar way.

\medskip\noindent
{\large\bf Acknowledgment}

This research was supported by the grant
from the National Natural Science Foundation of China 11571220.
Li is an honorary professor of the University of Hong Kong and the
Shanghai University. His research was also supported by USA NSF grant DMS
1331021, Simons Foundation Grant 351047.

\end{document}